\documentstyle[12pt]{amsart}

\newtheorem{thm}{Theorem}
\newtheorem{lem}[thm]{Lemma}
\newtheorem{cor}[thm]{Corollary}
\newtheorem{pr}[thm]{Proposition}

\newcommand{\N}{{\Bbb{N}}}
\newcommand{\R}{{\Bbb{R}}}
\newcommand{\al}{\alpha}
\newcommand{\om}{\omega}
\newcommand{\ep}{\epsilon}
\newcommand{\de}{\delta}
\newcommand{\sj}{\sum^m_{j=1}}
\newcommand{\sun}{\sum^n_{i=1}}
\newcommand{\su}{\sum^\infty_{i=1}}
\newcommand{\sk}{\sum^\infty_{k=1}}
\newcommand{\sik}{\sum^{i_k}_{i=i_{k-1}+1}\!}
\newcommand{\sil}{\sum^{i_l}_{i=i_{l-1}+1}\!}
\newcommand{\lam}{\lambda}
\newcommand{\co}{\mbox{co}}

\newcommand{\la}{\langle}
\newcommand{\ra}{\rangle}

\begin{document}

\title{Embedding Orlicz Sequence Spaces into $C(\alpha)$}
\author{Denny H. Leung}
\address{Department of Mathematics \\ National University of Singapore\\
         Singapore 0511}
\email{matlhh@@leonis.nus.sg}
\maketitle

\begin{abstract}
  Let $M$ be a non-degenerate Orlicz function such that there exist
$\ep > 0$ and $0 < s < 1$ with  $\su M(\ep s^i)/M(s^i) <
\infty$. It is shown that the Orlicz sequence space $h_M$ is
isomorphic to a subspace of $C(\om^\om)$.  It is also shown
that for any non-degenerate Orlicz function $M$, $h_M$ does not embed
into $C(\al)$ for any $\al < \om^\om$.
\end{abstract}

\section{Introduction}
Let $\al$ be an ordinal, then $\al+1$, the set of all ordinals preceding
or equal to $\al$, is a compact Hausdorff topological space under the
order topology.  Following common practice, we let $C(\al)$ denote the
space of all real-valued continuous functions defined on $\al+1$.  In
this note, we show that a certain class of Orlicz sequence spaces
embeds into $C(\om^\om)$.  It is also shown that no non-trivial Orlicz
sequence space embeds into $C(\al)$ for any $\al < \om^\om$.

We follow throughout standard Banach space terminology as may be found
in \cite{LT}.  An {\em Orlicz function}\/ is a real-valued
convex non-decreasing
function $M$ defined on $[0,\infty)$ such that $M(0) = 0$ and $M(t) >
0$ for some $t$.  It is
{\em non-degenerate}\/ if $M(t) > 0$ for all $t > 0$.  Given an Orlicz
function $M$, the Orlicz sequence space $h_M$ is the space of all real
sequences $(a_i)^\infty_{i=1}$ such that $\su M(c|a_i|) < \infty$ for
all $c < \infty$, endowed with the norm
\[ \|(a_i)\| = \inf\{\rho > 0 : \su M(|a_i|/\rho) \leq 1 \}.
\]
Two Orlicz functions $M$ and $N$ are said to be {\em equivalent}\/ if
there exist positive constants $K, k$ and $t_0$ such that 
\[ K^{-1}M(k^{-1}t) \leq N(t) \leq KM(kt) \]
whenever $0 \leq t \leq t_0$.
This holds if and only if $h_M = h_N$ as sets, and the identity map is
an isomorphism between the spaces \cite{LT}. Let $M$ be a continuously
differentiable Orlicz function such that $M'$ is a bijection from
$[0,\infty)$ onto itself.  Define $M^*(t) = \int^t_0(M')^{-1}(u) du$
for all $t \geq 0$.  Then $M^*$ is an Orlicz function.  Furthermore,
the following conditions are equivalent \cite{LT}:
\begin{enumerate}
\item $h_M$ does not contain an isomorph of $\ell^1$,
\item $h'_M$ is naturally isomorphic to $h_{M^*}$ via the duality
\[ \la(a_i), (b_i)\ra = \su\!a_ib_i \]
for all $(a_i) \in h_M$, and $(b_i) \in h_{M^*}$,
\item $M^*$ satisfies the $\Delta_2$-condition, i.e., $\limsup_{t\to
0}M^*(2t)/M^*(t) < \infty$.
\end{enumerate}

Equip the set of all real sequences with the coordinatewise partial
order. If $b = (b_i)$, let $|b| = (|b_i|)$. 
A set $S$ of such sequences is {\em solid}\/ if whenever $|y| \leq
|x|$ for some $x \in S$, then $y \in S$.
Finally, $c_{00}$ denotes the space of all finitely non-zero real
sequences, and the cardinality of a set $A$ is denoted by $|A|$.

%% section 1

\section{Embedding $h_M$ into $C(\om^\om)$}

Throughout this section, let $M$ be a fixed continuously
differentiable Orlicz function such that $M'$ is a bijection from 
$[0,\infty)$ onto itself, and let $t_n =
(M^*)^{-1}(1/n)$ for all $n \in \N$. 
It is also assumed that there
are $\ep > 0$ and $0 < s < 1$ such that $\su\! M(\ep s^i)/M(s^i)
< \infty$.  Note that if $s^i < t \leq s^{i-1}$, then
\[ \frac{M(\ep st)}{M(t)} \leq \frac{M(\ep s^i)}{M(s^i)} \to 0 \]
as $i \to \infty$.  Hence $\lim_{t\to 0}M(\ep st)/M(t) = 0$.  It
follows (see \cite{L}) that $h_M$ does not contain a copy of $\ell^1$.
By the remarks in  the Introduction, $M^*$ satisfies the
$\Delta_2$-condition, and $h_{M^*}$ is naturally isomorphic to $h'_M$. 
The main step in the embedding of $h_M$ into $C(\om^\om)$
lies in the solution of the optimization problem contained in the
following proposition. 

\begin{pr}\label{main}
For any $n \in \N$, let $S$ be the set of all non-negative,
non-increasing sequences $(b_i)^n_{i=1}$ such that $\sun\! M^*(b_i) =
1$.  Then there is a constant $C$, independent of $n$, such that,
taking $b_{n+1} = 0$, 
\[ \sun\frac{b_i - b_{i+1}}{t_i} \leq C \]
for all $(b_i)^n_{i=1} \in S$.
\end{pr}

\begin{pf}
We begin with some notation.
Define $f : S \to \R$ by
\[ f(b_1,\ldots, b_n) = \sun\frac{b_i - b_{i+1}}{t_i} .\]
Clearly, $f$ is a continuous function on the compact subset $S$ of
$\R^n$, hence it attains its maximum at a point, say $(a_i)^n_{i=1}$,
of $S$.  Let $0 = i_0 < i_1 < \cdots < i_m \leq n$ be such that 
\[ a_1 = \cdots = a_{i_1} > a_{i_1+1} = \cdots = a_{i_2} > \cdots >
a_{i_{m-1}+1} = \cdots = a_{i_m} > a_{{i_m}+1} = 0. \]

\begin{lem}\label{lambda}
For $2 \leq j \leq m$,
\[  \frac{1}{i_1t_{i_1}(M^*)'(a_{i_1})} = 
\frac{1}{(i_j-i_{j-1})(M^*)'(a_{i_j})}(\frac{1}{t_{i_{j}}} -
\frac{1}{t_{i_{j-1}}}). \]
Furthermore, labeling this common value by $\lam$, 
\[ f(a_1,\ldots,a_n) \leq \lam\sun M^*(2a_i) .\]
\end{lem}

\begin{pf}
For any $2 \leq j \leq m$, let $\de = \de(\ep)$ be defined in a
neighborhood of $0$ so that 
\begin{equation}\label{one}
\sum^{i_1}_{i=1}M^*(a_i+\ep) +
\sum^{i_j}_{i=i_{j-1}+1}\!\!M^*(a_i-\de(\ep))
+ \sum_{l\neq 1,j}\,\sum^{i_l}_{i=i_{l-1}+1}\!\! M^*(a_i) = 1 .
\end{equation}
Since $M^*$ is continuously differentiable, so is $\de$.  It is also
clear that $\de(\ep) \to 0$ as $\ep \to 0$.  Therefore, for $\ep$ in a
small enough neighborhood of $0$, the sequence
\[ x_\ep = (a_i)^n_{i=1} + \ep\sum^{i_1}_{i=1}e_i -
\de(\ep)\!\!\sum^{i_j}_{i=i_{j-1}+1}e_i \]
belongs to $S$, where $e_i$ is the element in $\R^n$ with $i$-th
coordinate $1$ and all other coordinates $0$.  By the maximality of
$f$ at $(a_i)$, we have
\[ \frac{d}{d\ep}f(x_\ep)|_{\ep=0} = 0 . \]
Hence
\begin{equation}\label{two}
\frac{1}{t_{i_1}} + \de'(0)(\frac{1}{t_{i_{j-1}}} - \frac{1}{t_{i_j}})
= 0 .
\end{equation}
Differentiating equation (\ref{one}) with respect to $\ep$, and
setting $\ep = 0$, we see that 
\[ \de'(0) = \frac{i_1(M^*)'(a_{i_1})}{(i_j-i_{j-1})(M^*)'(a_{i_j})} . \]
Combining this with equation (\ref{two}), we have 
\[  \frac{1}{i_1t_{i_1}(M^*)'(a_{i_1})} =
\frac{1}{(i_j-i_{j-1})(M^*)'(a_{i_j})}(\frac{1}{t_{i_{j}}} -
\frac{1}{t_{i_{j-1}}}) \]
for $2 \leq j \leq m$.  For convenience, let 
$\eta_1 = t_{i_1}$, and 
$\eta_j = 
(t^{-1}_{i_{j}} - t^{-1}_{i_{j-1}})^{-1}$ for $2 \leq j \leq m$.  Then
\begin{eqnarray*}
f(a_1,\ldots,a_n) & = & 
\sj\frac{{a_{i_j}}}{\eta_j} \\ 
& = & \lam\sum^m_{j=1}a_{i_j}(i_j - i_{j-1})(M^*)'(a_{i_j}) \\
& \leq & \lam\sum^m_{j=1}(i_j - i_{j-1})M^*(2a_{i_j}) \\
& = & \lam\sun M^*(2a_{i}),
\end{eqnarray*}
as claimed.
\end{pf}
Continuing with the proof of Proposition \ref{main}, we let $N(t) =
M^*(t)/t$  for all $t > 0$.  Then  $N$ is an increasing
function on $(0, \infty)$. Note that 
\[t_n = M^*(t_n)/N(t_n) = \frac{1}{nN(t_n)} \]
for all $n \in \N$. Hence $\eta^{-1}_1 =
i_1N(t_{i_1})$ and 
\begin{eqnarray}\label{three} 
\frac{1}{\eta_j} & = & \frac{1}{t_{i_{j}}} - \frac{1}{t_{i_{j-1}}}
        \nonumber \\ 
   & = & i_{j}N(t_{i_{j}}) - i_{j-1}N(t_{i_{j-1}}) \\
   & \leq & (i_{j} - i_{j-1})N(t_{i_{j}}) . \nonumber 
\end{eqnarray}
for $2 \leq j \leq m$. Now
\begin{eqnarray*}
1 & = & \sun M^*(a_i) \\
  & = & \sum^m_{j=1}(i_j - i_{j-1})M^*(a_{i_j}) \\
  & \leq & \sj (i_j - i_{j-1})a_{i_j}(M^*)'(a_{i_j}) \\
  & = & \sj\frac{a_{i_j}}{\lam\eta_j} \\
  & = & \sj\frac{1}{\lam\eta_j}M'(\frac{1}{\lam(i_j - i_{j-1})\eta_j})
               \hspace{2em} \mbox{since $(M^*)' = (M')^{-1}$} \\
  & \leq & \sj\frac{(i_j -
i_{j-1})N(t_{i_j})}{\lam}M'(\frac{N(t_{i_j})}{\lam}) 
               \hspace{2em} \mbox{by (\ref{three})}    \\
  & \leq & \sj(i_j - i_{j-1})M(\frac{2N(t_{i_j})}{\lam}) \\
  & \leq & \su M(\frac{2N(t_i)}{\lam}) .
\end{eqnarray*}
Let $K$ be the largest integer, possibly negative, such that  $N(t_1)
\leq s^K$.  By the assumption on $M$, $\lim_{t\to 0}(M')^{-1}(t) = 0$.
Hence $\lim_{t\to 0}N(t) = 0$.
For each $k \in \N$, let $l_k$ be the largest integer such
that $N(t_{l_k}) > s^{K+k}$.  Note that $s^{K+k} < N(t_{l_k}) =
(l_kt_{l_k})^{-1}$.  Hence $l_k < (s^{K+k}t_{l_k})^{-1}$.  Also,
$s^{K+k} < N(t_{l_k}) \leq (M^*)'(t_{l_k})$.  Therefore,
\[ t_{l_k} > M'(s^{K+k}) \geq M(s^{K+k})/s^{K+k} .\]
It follows that 
\begin{eqnarray*}
1 & \leq & \su M(\frac{2N(t_i)}{\lam}) \\
  & \leq & \sk l_kM(\frac{2}{\lam}s^{K+k-1}) \\
  & < & \sk \frac{1}{s^{K+k}t_{l_k}}M(\frac{2}{\lam}s^{K+k-1}) \\
  & < & \sk\frac{M(\frac{2}{\lam}s^{K+k-1})}{M(s^{K+k})} .
\end{eqnarray*}
From the assumption on $M$, one easily deduces that $\sk M(\ep
s^{K+k})/M(s^{K+k}) < \infty$.   
Consequently, if 
\[ \lam \geq \max(\frac{1}{s\ep}\sk\frac{M(\ep
s^{K+k})}{M(s^{K+k})},\frac{2}{s\ep})\, , \]
then
\begin{eqnarray*}
1 & < & \sk\frac{M(\frac{2}{\lam}s^{K+k-1})}{M(s^{K+k})} \\
  & = & \sk\frac{M(\frac{2}{\lam s\ep}\ep s^{K+k})}{M(s^{K+k})} \\
  & \leq & \frac{2}{\lam s\ep}\sk\frac{M(\ep s^{K+k})}{M(s^{K+k})}
\\
  & \leq & 1,
\end{eqnarray*}
which is impossible.  Hence $\lam$ is bounded by a constant which is
independent of $n$.  
Since $a_i \leq t_1$ for all $i$, we obtain from Lemma \ref{lambda}
that 
\[ \sun\frac{a_i - a_{i+1}}{t_i} \leq \lam\sun\!M^*(2a_i) \leq
\lam\sup_{0<t\leq t_1}\frac{M^*(2t)}{M^*(t)}. \]
As it has been observed that $M^*$ satisfies the $\Delta_2$-condition,
the proof is complete.
\end{pf}

We now proceed to show that $h_M$ embeds in $C(K)$ for some countable
compact set $K$.  

\begin{lem}\label{S}
For all $n \in \N$, 
\[ S \subseteq C\co(\{t_j\sum^j_{i=1}e_i : 1 \leq j \leq
n\}\cup\{0\}), \]
where the constant $C$ is as given in Proposition \ref{main}, and
$\co$ denotes the convex hull.
\end{lem}

\begin{pf}
Let $(b_1,\ldots,b_n) \in S$.  Taking $b_{n+1} = 0$, define
\[ c_i = \frac{b_i - b_{i+1}}{t_i} \]
for $1 \leq i \leq n$.  Then $\sun c_i \leq C$ by Proposition
\ref{main}, and 
\[ \sun b_ie_i = \sun c_i(t_i\sum^i_{j=1}e_j) .\]
The result follows immediately.
\end{pf}

For each $n \in \N$, let $K_n$ be the set of all infinite real
sequences $(b_i)$ such 
that $(|b_i|) = t_n\sum_{i\in A}e_i$ for some $A \subseteq \N$ of
cardinality $\leq n$  (the empty sum is defined to be the zero
sequence). 

\begin{lem}\label{norm}
Let $K = \bigcup^\infty_{n=1}\!K_n$, and let $U$ be the unit ball of
$h_{M^*}$.  Then $C^{-1}U$ is contained in the weak* closed convex hull
of $K$.
\end{lem}

\begin{pf}
For each $n \in \N$, it is easy to see that $\co(K_n)$ is solid.  Hence
so is $\co(K)$.  Let $(b_i)$ be a non-negative sequence in $c_{00}$
such that $\|(b_i)\|_{M^*} = 1$.  Then $\su M^*(b_i) = 1$.  Let $\pi$
be a permutation on $\N$ so that $(b_{\pi(i)})$ is non-increasing.  By
Lemma \ref{S}, 
\[ \sun b_{\pi(i)}e_i \in  C\co(\{t_j\sum^j_{i=1}e_i : 1 \leq j \leq
n\}\cup\{0\}) ,\]
where $n$ is the largest integer such that $b_{\pi(n)} > 0$.  Hence
\[ (b_{i}) \in  C\co(\{t_j\sum^j_{i=1}e_{\pi(i)} : 1 \leq j \leq
n\}\cup\{0\}) \subseteq C\co(K).\]
Since $\co(K)$ is solid, it follows that $c_{00} \cap U \subseteq
C\co(K)$.  The required result follows immediately.
\end{pf}

\begin{pr}\label{embed}
The set $K$ is countable and compact under the weak* topology.
Moreover, $h_M$ is isomorphic to a subspace of $C(K)$.
\end{pr}

\begin{pf}
The countability of $K$ is obvious.
From the definition, it is easily seen that $K$ is contained in the
unit ball $U$ of $h_{M^*}$.  To show that $K$ is weak* compact, it
suffices to show that it is weak* closed.  We claim that in fact 
$\overline{K}^{w^*}= \bigcup^\infty_{n=1}\overline{K_n}^{w^*}$. 
Since it is clear that each $K_n$ is
weak* closed, the result follows.
Suppose now that  there exists $x \in
\overline{K}^{w^*}\backslash\bigcup^\infty_{n=1}\overline{K_n}^{w^*}$.  
Fix $i \in \N$.  For any $\ep > 0$ and $n \in \N$,
there exists $y \in K \backslash \bigcup^n_{j=1}K_j$ such that 
$|\la e_i, x - y\ra| < \ep$.  Since $|\la e_i, y\ra| \leq t_{n+1}$, we
see that $|\la e_i, x\ra| \leq \ep + t_{n+1}$.  As $\ep > 0$ is
arbitrary, and $t_n \to 0$ as $n \to \infty$, we see that $x = 0$,
which is impossible.

By virtue of Lemma \ref{norm}, the evaluation map $T : h_M \to C(K)$,
$(Tb)(x) = \la b, x\ra$ for all $x \in K$, is an isomorphic embedding.
\end{pf}

Recall that for a compact Hausdorff space $X$, the {\em derived set}\/
$X^{(1)}$ is the set of all limit points of $X$. Define a transfinite
inductive sequence of sets as follows: $X^{(0)} = X$, if $X^{(\al)}$
has been defined for all ordinals $\al < \beta$, let $X^{(\beta)}$ be
the derived set of $X^{(\gamma)}$ if $\beta = \gamma + 1$; otherwise,
let $X^{(\beta)} = \bigcap_{\al<\beta}X^{(\al)}$.

\begin{pr}\label{om}
The topological space $K$ is homeomorphic to $\om^\om +1$.
\end{pr}

\begin{pf}
By a result of Semadeni, (cf. \cite[Corollary 5.2]{La}), it suffices
to show that $K^{(\om)}$ is a singleton.  We claim that in fact it
consists of $0$ only. From the proof of Proposition \ref{embed}, we
see that $K^{(1)} = \bigcup^\infty_{n=1}K^{(1)}_n$.  It is easy to see
that 
\[ K^{(1)}_n = \{(b_i) : (|b_i|) = t_n\sum_{i\in A}e_i, |A| \leq
n-1\}. \] 
Repeating the argument, we obtain that $K^{(m)} =
\bigcup^\infty_{n=m}K^{(m)}_n$, and that 
\[K^{(m)}_n = \{(b_i) :
(|b_i|) = t_n\sum_{i\in A}e_i, |A| \leq n-m\} \]
for all $m \leq n \in \N$.
Consequently,  $K^{(\om)} = \bigcap^\infty_{m=1}K^{(m)} = \{0\}$, as
required. 
\end{pf}

Finally, we arrive at the following result.

\begin{thm}
Let $N$ be a non-degenerate Orlicz function such that there exist $\ep
> 0$ and $0 < s < 1$ such that $\su N(\ep s^i)/N(s^i) < \infty$.  Then
$h_N$ is isomorphic to a subspace of $C(\om^\om)$.
\end{thm}

\begin{pf}
It is well known that $N$ is equivalent to a continuously
differentiable Orlicz function $M_1$ (see \cite{LT}).  
Since $h_N$ does not contain a copy of $\ell^1$, neither does
$h_{M_1}$. Hence $M'_1(0) =0$. By the non-degeneracy, $M'_1(t) > 0$
for all $t > 0$.  Let $M(t) = \int^t_0(1+u)M'_1(u) du$ for $t \geq 0$.
It is easy to see that $M$ satisfies all the conditions stated at the
beginning of the section.  Furthermore, $M$ is equivalent to $M_1$,
and thus $N$, because
\[ M_1(t) \leq M(t) \leq 2M_1(t) \]
for $0 \leq t \leq 1$.
By Propositions
\ref{embed} and \ref{om}, $h_M$ embeds in $C(\om^\om)$.  Since $h_N$
and $h_M$ are isomorphic, the result follows.
\end{pf}

%section 2

\section{A non-embedding result}
  
In this section, we show that if $M$ is a non-degenerate Orlicz
function, then $h_M$ does not embed into $C(\al)$ for any $\al <
\om^\om$. We begin with the following well known lemma.

\begin{lem}\label{easy}
A Banach space $E$ embeds into $c_0$
if and only if there is a weak* null sequence $(x'_i)$ in
$E'$, and a constant $C > 0$ such that
\[ \|x\| < C\sup_i|\la x, x'_i\ra| \]
for all $x \in E$.
\end{lem}

\begin{thm}\label{non}
Let $E$ be a Banach space with a shrinking basis $(e_i)$.  Then $E$
embeds into $c_0$ if and only if there is a strictly increasing
sequence $(i_k)^\infty_{k=1}$ in $\N$ such that 
$\su a_ie_i$ converges in $E$ whenever $\|\sik a_ie_i\| \to 0$ as $k
\to \infty$, where $i_0 = 0$. 
\end{thm}

\begin{pf}
Let the sequence $(i_k)$ have the property stated.  For all $k \in
\N$, let $E_k = [e_i]^{i_k}_{i=i_{k-1}+1}$.  Then $E$ is isomorphic to
$(\sum\!\oplus E_k)_{c_0}$ by the Open Mapping Theorem.  
For each $k$, choose a finite subset $W_k$
of the unit ball of $E'_k$ such that 
\[ \|x\| \leq 2\sup_{x'\in W_k}|\la x, x'\ra| \]
for all $x \in E_k$.  Using the obvious identification of $E'_k$ with
a subset in $E'$, we can arrange $\bigcup W_k$ into a weak* null
sequence $(x'_i)$ satisfying the condition in Lemma \ref{easy}. Hence
$E$ embeds into $c_0$, as required.

Conversely, assume that $E$ embeds into $c_0$.  Choose a sequence
$(x'_i)$ and a constant $C$ as in Lemma \ref{easy}.  Let $(e'_i)$ be
the sequence biorthogonal to $(e_i)$.  Since $(e_i)$ is shrinking,
$(e'_i)$ is a basis of $E'$. Let $(P_i)$ be the basis projections with
respect to $(e_i)$.  
As $(x'_i)$ is weak* null, we see that $\lim_{i\to\infty}P'_jx'_i
= 0$ in norm for any fixed $j$. Applying an easy perturbation
argument, we may assume the existence of two non-decreasing sequences
of natural numbers, $(m_k)$ and $(n_k)$ such that $m_k \leq n_k$ for
all $k \in \N$, $\lim_{k\to\infty}m_k = \lim_{k\to\infty}n_k =
\infty$, and $x'_k \in [e'_i]^{n_k}_{i=m_k}$ for all $k$.
Choose a strictly increasing sequence $(k_l)$ in $\N$ inductively, so
that $k_1 = 1$, and $n_{k_l} < m_{k_{l+1}}$ for all $l \in \N$.
For each $l \in \N$, let 
\[ V_l = \{(P_{n_{k_l}}-P_{n_{k_{l-1}}})'x'_i : k_{l-1} < i < k_{l+1}\} , \]
where $n_0 = k_0 = 0$.  Clearly $V_l \subseteq
[e'_i]^{n_{k_l}}_{i=n_{k_{l-1}}+1}$. Also, each $x'_i$ can be written as a
sum of at most two elements from $\bigcup^\infty_{l=1}V_l$.  Hence
\[ \|x\| < \frac{C}{2}\sup\{|\la x, v'\ra| : v' \in
\bigcup^\infty_{l=1}V_l \} . \]
Finally, let $i_l = n_{k_l}$, $l \in \N$.  
Then for any 
$(a_i) \in c_{00}$, 
\begin{eqnarray*}
\|\su a_ie_i\| & \leq & 
\frac{C}{2}\sup_l\sup_{v'\in V_l}|\la \su a_ie_i, v'\ra|
\\
 & \leq & \tilde{C}\sup_l\|\!\sil a_ie_i\|
\end{eqnarray*}
for some constant $\tilde{C}$. The desired conclusion follows easily.
\end{pf}

\begin{cor}
Let $M$ be a non-degenerate Orlicz function, then $h_M$ does not embed
into $C(\al)$ for any ordinal $\al < \om^\om$.
\end{cor}

\begin{pf}
Certainly, it is only necessary to consider infinite ordinals.  It is
well known that $C(\al)$ is isomorphic to $c_0$ for any infinite
ordinal $\al < \om^\om$. Thus if $h_M$ embeds into $C(\al)$
for some such $\al$, it also embeds into $c_0$.  In particular, the
coordinate unit vectors $(e_i)$ form a shrinking basis of $h_M$.  We
then obtain a sequence $(i_k)$ having the properties enunciated in
Theorem \ref{non}. For convenience, we assume the normalization
condition $M(1) = 1$.  Choose a sequence of finite subsets $(A_j)$ of
$\N$ such that $\max A_j < \min A_{j+1}$  and
$|A_j|M(1/j) \geq 1$ for all $j \in \N$.  Note that the second
condition requires the non-degeneracy of $M$.  The formal sum
$\sum^\infty_{j=1}\sum_{k\in A_j}e_{i_k}/j$ satisfies the condition in
Theorem \ref{non} and hence converges in $h_M$.  But this is
impossible since $\sum^\infty_{j=1}|A_j|M(1/j) = \infty$.
\end{pf}  

We close with the following natural question concerning the results of
this paper.\\

\noindent{\em Problem\ }
Identify all Orlicz functions $M$
such that $h_M$ embeds into $C(\al)$ for some ordinal $\al$.  \\

Let us remark that as  a consequence of the results in \cite{L}, there
is an Orlicz function $M$ such that $h_M$ is $c_0$-saturated, but does
not embed into $C(\al)$ for any ordinal $\al$.

%bibliography

\bibliographystyle{standard}
\bibliography{tref}

\end{document}